\documentclass[conference, 10pt]{IEEEtran}

\usepackage{cite}
\usepackage[font=footnotesize]{subfig}
\usepackage{url}
\usepackage{graphicx}
\usepackage{color}
\usepackage{placeins}
\usepackage{float}
\usepackage{tabularx,colortbl}
\usepackage{ifthen}

\usepackage{amsmath,amssymb,amsfonts}
\usepackage{algorithmic}
\usepackage{textcomp}
\usepackage{siunitx}

\newcommand{\diff}{\,\mathrm{d}}

\newcommand{\nh}{\vv{\hat{n}_{{r}}}}

\newcommand{\Tc}{ { \mathcal{T}} }

\newcommand{\Kc}{ { \mathcal{K}} }
\newcommand{\Ic}{ { \mathcal{I}} }
\newcommand{\Pc}{ { \mathcal{P}} }
\newcommand{\Rc}{ { \mathcal{R}} }

\newcommand{\vv}[1]{\ensuremath{\boldsymbol{#1}}}                   

\newcommand{\vet}[1]{{ \mathbf{#1}}}                      

\newcommand{\mat}[1]{\mathbf{{#1}}}         
\newcommand{\matb}[1]{\mathbb{{#1}}}         

\makeatletter

\newcounter{author}
\renewcommand{\author}[2][]{
   \stepcounter{author}
   \@namedef{author@\theauthor}{#2}
   \@namedef{authorlabel@\theauthor}{#1}
}

\newcounter{address}
\newcommand{\address}[2][]{
   \stepcounter{address}
   \@namedef{address@\theaddress}{#2}
   \@namedef{addresslabel@\theaddress}{#1}
}

\newcommand{\alsep}{and}

\def\newmaketitle{\par%
  \begingroup%
  \normalfont%
  \def\thefootnote{}
  \def\footnotemark{}
  \let\@makefnmark\relax
  \footnotesize
  \footnotesep 0.7\baselineskip
  \normalsize%
  \twocolumn[\thenewmaketitle\@IEEEaftertitletext]%
  \if@IEEEusingpubid
     \enlargethispage{-\@IEEEpubidpullup}%
  \fi
  \endgroup
  \setcounter{footnote}{0}\let\maketitle\relax\let\@maketitle\relax
  \gdef\@thanks{}%
  \let\thanks\relax}

\def\thenewmaketitle{
  \newpage
  \begin{center}%
    \vskip0.2em{\Huge\@IEEEcompsoconly{\sffamily}\@IEEEcompsocconfonly{\normalfont\normalsize\vskip 2\@IEEEnormalsizeunitybaselineskip
   \bfseries\large}\@title\par}\vskip1.0em\par%
    \vspace{1ex}
    \newcounter{c@author}
    \newcounter{c@tmp}
    \ifthenelse{\value{author}=2}{%
      \newcommand{\liand}{ and }}{%
      \newcommand{\liand}{, and }}
    \ifthenelse{\value{address}<2}{%
      \@nameuse{author@1}%
      \stepcounter{c@author}%
      \whiledo{\value{c@author}<\value{author}}{%
        \setcounter{c@tmp}{\value{author}}%
        \addtocounter{c@tmp}{-\value{c@author}}%
        \ifthenelse{\value{c@tmp}=1}{%
          \renewcommand{\alsep}{\liand}}{\renewcommand{\alsep}{, }}%
        \stepcounter{c@author}\alsep \@nameuse{author@\thec@author}}\\%
    }
    {
      \@nameuse{author@1}${}^{(\ref{\@nameuse{authorlabel@1}})}$%
      \stepcounter{c@author}%
      \whiledo{\value{c@author}<\value{author}}{%
      \setcounter{c@tmp}{\value{author}}%
      \addtocounter{c@tmp}{-\value{c@author}}%
      \ifthenelse{\value{c@tmp}=1}{%
        \renewcommand{\alsep}{\liand}}{\renewcommand{\alsep}{, }}%
      \stepcounter{c@author}\alsep \@nameuse{author@\thec@author}%
        ${}^{(\ref{\@nameuse{authorlabel@\thec@author}})}$%
      }
    }
    \vspace{0.2ex}

    \ifthenelse{\value{address}>0}{%
      \ifthenelse{\value{address}=1}{
        {\@nameuse{address@1}}
      }
      {
        \newcounter{c@address}

        \begin{center}
        \whiledo{\value{c@address}<\value{address}}
        {
          \refstepcounter{c@address}
            ${}^{(\thec@address)}$\,%
              \label{\@nameuse{addresslabel@\thec@address}}%
              \@nameuse{address@\thec@address}\\ %
        }
        \end{center}
      } 
    }
    {
      \relax
    }
  \end{center}
}

\makeatother

\title{On a Frequency-Stabilized Single Current Inverse Source Formulation}

\author[org1]{Paolo Ricci}
\author[org2]{Adrien Merlini}
\author[org1]{Francesco P. Andriulli}

\address[org1]{Politecnico di Torino, Turin, Italy}
\address[org2]{IMT Atlantique, Brest, France}

\begin{document}

\newmaketitle

\begin{abstract}
Several strategies are available for solving the inverse source problem in electromagnetics. Among them, many have been focusing in retrieving Love currents by solving, after regularization, for Love's electric and magnetic currents. In this work we present a dual-element discretization, analysis, and stabilization of an inverse source formulation providing Love data by solving for only one current. This results in substantial savings and allows for an effective quasi-Helmholtz projector stabilization of the resulting operator. Theoretical considerations are complemented by numerical tests showing effectiveness and efficiency of the newly proposed method.
\end{abstract}

\section{Introduction}

Inverse source strategies based on external measurements are widely used for several applications including antenna diagnostics and  near-field (NF) to far-field (FF) transformations.
These methods are challenging because of the ill-posedness of the radiation operator.
The different approaches present in literature can be divided into two categories: methods that consider both electric and magnetic currents and deal with the null-space introduced by enforcing an additional constraint to the solution---often the Love (or zero-field inside) condition \cite{quijano2010field}---and strategies that partially reduce the space of solutions by restricting the unknowns to a single type of current \cite{quijano2010field,kornprobst2021accuracy}.
In both cases, the Love conditions can be enforced. While the Love constraints have no impact on the reconstruction of the FF, their enforcement can be useful for diagnostic purposes and can be performed either explicitly or implicitly by leveraging an approximation of the Steklov-Poincaré operator \cite{kornprobst2021accuracy}.

In this work we propose a scheme falling in the second category that targets the reconstruction of only one of the two Love currents. Differently from previous works, however, we propose a single source approach that does not rely on any approximate impedance condition, but directly targets a stable discretization of the Steklov-Poincare operator that leverages on dual elements. This has two advantages: in addition to providing Love currents without scattering approximations, it enables the stabilization of the formulation at lower frequencies. A second contribution of this work is, in fact, a low-frequency stabilization of the Steklov-Poincaré operator based on quasi-Helmholtz projectors. The effectiveness of the new schemes is validated both by theoretical considerations and by numerical results that show the relevance of the methods in real case scenarios. 

\section{Notation and Background}

Let $\Gamma$ be a closed, simply connected surface surrounding an electric source.
By means of the equivalence theorem, equivalent electric and magnetic surface currents densities $\vv{J}$ and $\vv{M}$ can be obtained on $\Gamma$ that radiate the same external field as the original source.
We define the operators
\begin{align}
\Tc(\vv{J})(\vv{r}) &= k\, \mathcal{T}_s(\vv{J})(\vv{r})+ k^{-1}\, \mathcal{T}_h(\vv{J})(\vv{r})\,,\\
\Kc(\vv{J})(\vv{r}) &=  - \nh \times p.v. \int_\Gamma
\nabla \times \frac{e^{i k \left|\vv{r}-\vv{r}'\right|}}{4 \pi \left|\vv{r}-\vv{r}'\right|} \vv{J}(\vv{r}') \diff s'
\end{align}
with $\mathcal{T}_s(\vv{J})(\vv{r})=i\nh\times\int_{\Gamma}\frac{e^{i k\left|\vv{r}-\vv{r}'\right|}}{4\pi\left|\vv{r}-\vv{r}'\right|}{\vv{J}(\vv{r}')}\diff \vv{r}'$, $\mathcal{T}_h(\vv{J})(\vv{r})=i\nh\times\nabla\int_{\Gamma}\frac{e^{ik\left|\vv{r}-\vv{r}'\right|}}{4\pi\left|\vv{r}-\vv{r}'\right|}{\nabla_s\cdot\vv{J}(\vv{r}')}\diff \vv{r}'$,
and where $k$ is the wavenumber and $\nh$ is outside pointing normal vector at $\vv{r} \in \Gamma$.
Let $\Pc^+$ and $\Pc^-$ be respectively the external and the internal Calder\'on projectors
that satisfy the property $\Pc^-\Pc^+=\vv{0}$, from which we obtain that
\begin{equation}
    \Pc^-
    \begin{bmatrix}
    -\vv{M}^e\\
    \eta_0 \vv{J}^e
    \end{bmatrix}
    = 
    \begin{bmatrix}
    \frac{\Ic}{2}+\Kc&-\Tc\\
    \Tc&\frac{\Ic}{2}+\Kc
    \end{bmatrix}
    \begin{bmatrix}
    -\vv{M}^e\\
    \eta_0 \vv{J}^e
    \end{bmatrix}
    = \vv{0}
    \label{eq:zeroprop}
\end{equation}
where $\vv{M}^e$ and $\vv{J}^e$ are Love magnetic and electric currents on $\Gamma$ (radiating zero-fields inside and the original field outside). Given a second measurement surface $\Gamma_m$ external to $\Gamma$, we denote by
\begin{equation}
    \Rc=\begin{bmatrix}
-\Kc_{\vv{r}}&\Tc_{\vv{r}}\\
-\Tc_{\vv{r}} &-\Kc_{\vv{r}}
\end{bmatrix}
    \label{eq:radmat}
\end{equation}
the tangential projection of the radiation operator that defines the relation between the currents on $\Gamma$ and the tangential electric and magnetic fields ($\nh\times\vv{E},\,\nh\times\vv{H}$) on the measurement surface $\Gamma_m$, with $\Tc_{\vv{r}}=\Tc$. The operator $\Kc_{\vv{r}}$ is similar to $\Kc$ without the principal value. The corresponding linear system is
\begin{equation}\label{system}
\Rc
\begin{bmatrix}
-\vv{M} \\
\eta_0\vv{J}
\end{bmatrix}
=
\begin{bmatrix}
\nh\times\vv{E} \\
\eta_0\nh\times\vv{H}
\end{bmatrix}\,.
\end{equation}

\section{A Stable Discretization of the Steklov-Poincaré Operator}

In a general inverse source setting, current sources can be recovered by pseudo-inverting (or iteratively pseudo-inverting) the discretized counterpart of $\Rc$ in \eqref{system} and  to which additional equations are added to enforce the Love condition \cite{quijano2010field}. In this work, instead, we consider a single row of \eqref{eq:radmat} (e.g. the first row, related to the electric field) and we formally solve for one of the currents which yields the equation
\begin{equation}
\left( -\Kc_{\vv{r}} - \Tc_{\vv{r}}\left(\frac{\Ic}{2} + \Kc\right)^{-1}\Tc \right) \left(-\vv{M}\right) = \nh\times\vv{E}\,,
\label{eq:stprad}
\end{equation}
where $\left(\frac{\Ic}{2} + \Kc\right)^{-1}\Tc$ is the Steklov-Poincar\'e operator.
Equation \eqref{eq:stprad} can, upon discretization, be pseudo-inverted using a Moore-Penrose pseudo-inverse to obtain the magnetic Love current.
In practice, however, the stable discretization of \eqref{eq:stprad} cannot be achieved by only using the well known Rao-Wilton-Glisson (RWG) functions. One possibility is to rely on an approximation of the Steklov-Poincar\'e operator as is done in \cite{kornprobst2021accuracy}. Here we will follow a different approach based on dual elements which, in addition to providing a consistent discretization of the operator, will also allow for its stabilization with quasi-Helmholtz projectors as will be discussed in the next section.
After approximating the geometries of $\Gamma$ and of $\Gamma_m$ with meshes of triangular elements with $N_e$ and $N_m$ edges respectively,
the discretization we propose is
\begin{equation}
\label{eq:discrete}
   \left(-\mat{K}_{\vv{r}} - \matb{T}_{\vv{r}} \left( \matb{G}/2 + \matb{K} \right)^{-1} \mat{T} \right) \left(-\vet{m}\right) =  \vet{e}
\end{equation}
where
$\left[ \vet{e} \right]_m = \left<\nh\times\vv{g}_m, \nh\times\vv{E}\right>_{\Gamma_m}$,
$\vv{M}(\vv{r})\approx \sum_{n=1}^{N_e} \left[ \vv{m} \right]_n \vv{f}_n(\vv{r})$,
$\vv{f}(\vv{r})$ are the RWG basis functions (defined without edge normalization),
$\vv{g}(\vv{r})$ are the Buffa-Christiansen basis functions (see \cite{andriulli2008multiplicative} for their definition),
and where
$[\mat{K}_{\vv{r}}]_{ij}=\left< \nh\times \vv{g}_i, \Kc_{\vv{r}} \vv{f}_j \right>_{\Gamma_m}$,
$[\matb{T}_{\vv{r}}]_{ij}=\left< \nh\times \vv{g}_i, \Tc \vv{g}_j\right>_{\Gamma_m}$,
$[\matb{G}]_{ij}=\left< \nh\times \vv{f}_i, \vv{g}_j\right>_\Gamma$,
$[\matb{K}]_{ij}=\left< \nh\times \vv{f}_i, \Kc \vv{g}_j \right>_\Gamma$, and
$[\mat{T}]_{ij}=\left< \nh\times \vv{f}_i, \Tc \vv{f}_j\right>_\Gamma$,
with $\left<\vv{a},\vv{b}\right>_\Gamma = \int_\Gamma\vv{a}(\vv{r})\cdot\vv{b}(\vv{r}) \diff s$.
After solving \eqref{eq:discrete} and recovering $\vet{m}$, the other Love current can be obtained by applying the Steklov-Poincar\'e operator once more.

\section{Low-Frequency Stabilization}

The operator in \eqref{eq:discrete} (as well as the original system \eqref{system}) is affected by low frequency breakdown, i.e. its condition number increases as the frequency decreases, leading to the insurgence of numerical errors.
The stabilization we propose is based on Quasi-Helmoltz projectors that are
defined as $\mat{P}^{\Sigma} = \mat{\Sigma} (\mat{\Sigma}^{\text{T}} \mat{\Sigma})^{+}\mat{\Sigma}^{\text{T}}$, $\mat{P}^{\Lambda H }=\mat{I}-\mat{P}^{\Sigma}$, $\matb{P}^{\Lambda} = \mat{\Lambda} (\mat{\Lambda}^{\text{T}} \mat{\Lambda})^{+}\mat{\Lambda}^{\text{T}}$, and  $\matb{P}^{\Sigma H}=\mat{I}-\matb{P}^{\Lambda}$, where $\mat{\Lambda}$ and $\mat{\Sigma}$ are respectively the loop- and the star-decomposition matrices defined in \cite{andriulli2012well}, and  $^+$  denotes the Moore-Penrose pseudo-inverse.
The stabilized equation reads
\begin{equation}
    \matb{M}
    \left(-\mat{K}_{\vv{r}} - (\matb{T}_{\vv{r}} \left( \matb{G}/2 + \matb{K} \right)^{-1} \mat{T} \right)
    \mat{M} \vet{x} = \matb{M} \vet{e}
    \label{eq:lfbreak}
\end{equation}
where $\mat{M}\vet{x} = -\vet{m}$, $\mat{M}=1/\sqrt{k}\mat{P}^{\Lambda H} + i\sqrt{k}\mat{P}^\Sigma$ and $\matb{M}$ is the same as $\mat{M}$ with the roles of $\mat{\Lambda}$ and $\mat{\Sigma}$ exchanged.
We have to omit the technical details for the sake of conciseness, but we can rigorously prove the stabilization by leveraging  the  cancellations of $\Tc_h$ on the solenoidal subspaces, the mutual orthogonality of loop and star coefficients, and finally on our proof of the property  $\lim_{k\to0}\| \matb{P}^\Lambda \left( \matb{G}/2 + \matb{K} \right)^{-1} \mat{P}^\Sigma\|= 0$.

\section{Numerical Results}

\begin{figure}
    \centering
    \includegraphics[width=.399\textwidth]{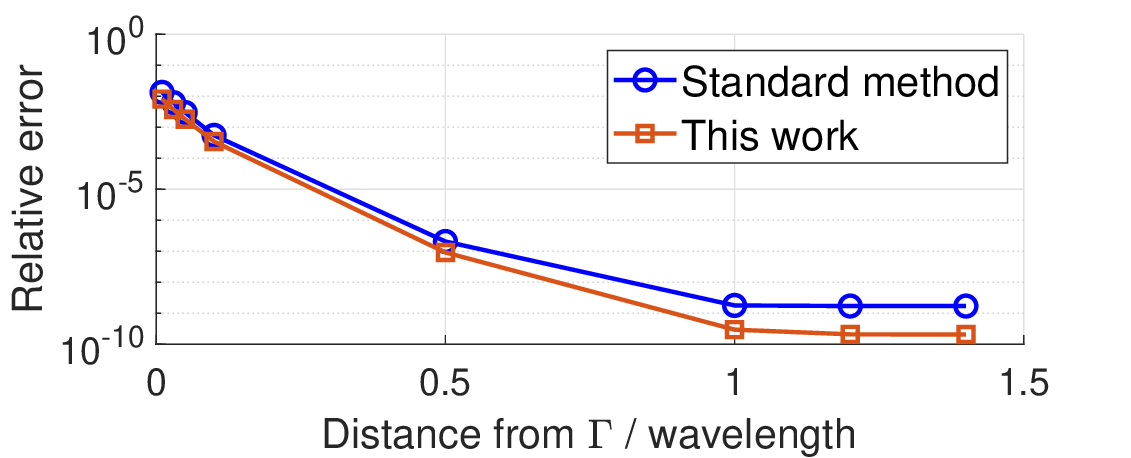}
    \caption{Relative error of the reconstruction of the fields over the distance from the equivalent surface in terms of wavelength.}
    \label{fig:fields}
\end{figure}

\begin{figure}
    \centering
    \includegraphics[width=.399\textwidth]{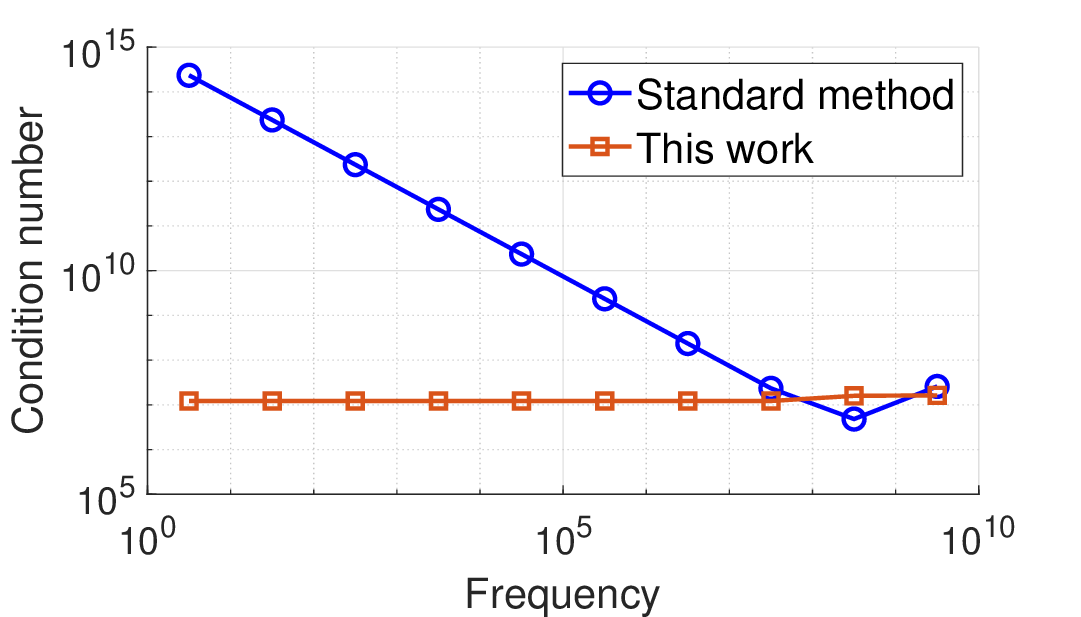}
    \caption{Condition number of the stabilized operator in function of the frequency.}
    \label{fig:sv}
    \vspace{-.5cm}
\end{figure}

The effectiveness of the approach we proposed has been verified by imaging an electric dipole oscillating at \SI{3.16}{\giga\hertz} located in the center of a spherical equivalent surface $\Gamma$ with a radius of \SI{4}{\centi\meter}.
First, the equivalent currents have been reconstructed from the fields scattered by the dipole at a distance of one wavelength from $\Gamma$; these currents have then been radiated on several concentric spherical surfaces with increasing radii. The error of the reconstructed fields with respect to the ones obtained from the dipole is depicted in Fig.~\ref{fig:fields} and clearly show that the Steklov-Poincaré-based approach achieves similar results in the reconstruction of the fields with respect to the popular standard formulation---based on the solution of \eqref{system} with Love constraints---despite being involving a significantly smaller linear system.
Furthermore we verified that the condition number of the stabilized operator evaluated at the pseudo-inversion threshold does not increase when the frequency decreases (Fig.~\ref{fig:sv}), showing that the new method is free from low frequency conditioning breakdown.

\section*{ACKNOWLEDGEMENT}
The work of this paper has received funding from the EU H2020 research and innovation programme under the Marie Skłodowska-Curie grant agreement n° 955476 (project COMPETE) and from the European Research Council (ERC) under the European Union’s Horizon 2020 research and innovation programme (grant agreement No 724846, project 321).



\bibliographystyle{IEEEtran}

\begin{thebibliography}{1}

\bibitem{quijano2010field}
J. L. A. Quijano and G. Vecchi, “Field and source equivalence in source
reconstruction on 3d surfaces,” \emph{Progress In Electromagnetics Research},
vol. 103, pp. 67–100, 2010.
\bibitem{kornprobst2021accuracy}
J. Kornprobst, J. Knapp, R. A. Mauermayer, O. Neitz, A. Paulus, and
T. F. Eibert, “Accuracy and conditioning of surface-source based near-
field to far-field transformations,” \emph{IEEE Transactions on Antennas and
Propagation}, vol. 69, no. 8, pp. 4894–4908, 2021.
\bibitem{andriulli2008multiplicative}
F. P. Andriulli, K. Cools, H. Bagci, F. Olyslager, A. Buffa, S. Christiansen,
and E. Michielssen, “A multiplicative calderon preconditioner for the
electric field integral equation,” \emph{IEEE Transactions on Antennas and
Propagation}, vol. 56, no. 8, pp. 2398–2412, 2008.
\bibitem{andriulli2012well}
F. P. Andriulli, K. Cools, I. Bogaert, and E. Michielssen, “On a well-
conditioned electric field integral operator for multiply connected geome-
tries,” \emph{IEEE transactions on antennas and propagation}, vol. 61, no. 4, pp.
2077–2087, 2012.
\end{thebibliography}

\end{document}